\documentclass[12pt]{amsart}
\usepackage{times,amscd,amssymb,fullpage,epsfig}
\usepackage[latin1]{inputenc}
\usepackage{color}
\definecolor{0}{rgb}{0.0,0.0,0.0}
\definecolor{2}{rgb}{0.5,0.0,0.1}
\definecolor{3}{rgb}{0.0,0.4,0.1}

\newtheorem{thm}{Theorem}

\newtheorem{lem}{Lemma}
\newtheorem{cor}{Corollary}

\theoremstyle{remark}

\newtheorem{ex}{Example}

\theoremstyle{definition}
\newtheorem{defn}{Definition}

\newcommand{\Aut}{\operatorname{Aut}}

\newcommand{\R}{\mathbb{ R}}

\newcommand{\Z}{\mathbb{ Z}}

\newcommand{\map}{\rightarrow}

\renewcommand{\geq}{\geqslant}
\renewcommand{\le}{\leqslant}
\renewcommand{\leq}{\leqslant}

\newcommand{\NN}{\mathbb N}

\newcommand{\ZZ}{\mathbb Z}

\newcommand{\Stab}{\mathop{\textrm{Stab}}\nolimits}

\title{The classification of football patterns}
\author{V.~Braungardt}
\address{Mathematisches Institut, Ludwig-Maximilians-Universit\"at M\"unchen,
Theresienstr.~39, 80333 M\"unchen, Germany}
\email{Volker.Braungardt@mathematik.uni-muenchen.de}
\author{D.~Kotschick}
\address{Mathematisches Institut, 
Ludwig-Maxi\-mi\-lians-Universit\"at M\"unchen,
Theresienstr.~39, 80333 M\"unchen, Germany}
\email{dieter@member.ams.org}
\date{April 5, 2006; MSC 2000 classification: primary 52B05, 57M12, 57M15; secondary 05C25, 20E08}

\begin{document}

\begin{abstract}
We prove that every spherical football is a branched cover, branched only in the vertices, of the standard football made up 
of $12$ pentagons and $20$ hexagons. We also give examples showing that the corresponding result is not true for 
footballs of higher genera. Moreover, we classify the possible pairs $(k,l)$ for which football patterns on the sphere exist
satisfying a natural generalisation of the usual incidence relation between pentagons and hexagons to $k$-gons and $l$-gons. 
\end{abstract}

\maketitle

\begin{center}
\epsfig{file=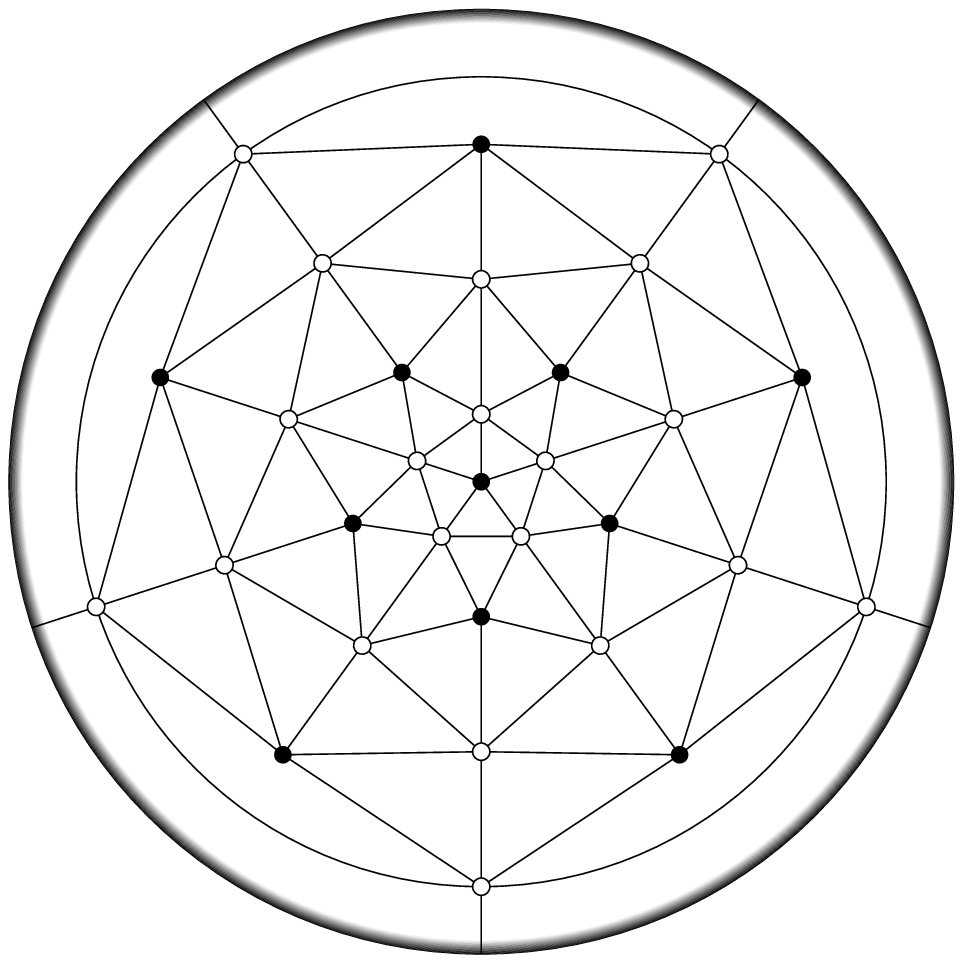}
\end{center}

\vfill\newpage

\section*{Introduction}

A football pattern\footnote{We use English terminology. American readers might want to call our football patterns 
``soccer ball patterns''.} is a graph embedded in the two-sphere in such a way that all faces are pentagons and hexagons, 
satisfying the conditions that 
the edges of each pentagon meet only edges of hexagons, and that the edges of each hexagon alternately meet 
edges of pentagons and of hexagons. If one requires that there are exactly three faces meeting at each vertex,
then Euler's formula implies that the pattern consists of $12$ pentagons and $20$ hexagons. Moreover, in this
case the combinatorics of the pattern is uniquely determined. This pattern, which we shall refer to as the
standard football, has a particularly symmetric realisation with all polygons regular, which can be thought of 
as a truncated icosahedron.

If one drops the requirement that there are exactly three polygons meeting at each vertex, then one can exhibit
infinitely many distinct football patterns by lifting to branched covers of the standard football branched only in
vertices of the pattern. In the first part of this paper we shall prove that these are the only football patterns on
the two-sphere. For the proof we consider the dual graph of a football pattern as a coloured ribbon graph. 
The dual of the standard football with $12$ pentagons and $20$ hexagons is shown in the picture on the 
first page of this paper. To make the picture symmetric, a black vertex is spread out at infinity and is depicted by
the shaded ring around the rest of the graph.

Of course we could consider the dual graph as a map, but in order to make clear the distinction between football 
patterns and their duals, we will always use the language of maps, cf.~\cite{Maps}, for the football patterns themselves, 
and the language of ribbon graphs, cf.~\cite{RibbonGraphs}, for their duals.

All the football graphs, that is ribbon graphs dual to football patterns, have the same universal covering, which is a 
certain tree $T$. We shall determine the automorphism group of $T$ and prove that every cofinite subgroup of the 
automorphism group whose quotient gives rise to the dual of a spherical pattern is a subgroup of the group giving rise 
to the standard football. We also consider football patterns on surfaces of higher genera and show that the classification 
theorem does not hold for them; in other words, not all of them can be obtained by taking branched covers of the standard 
football. However, football patterns  on surfaces of higher genera always admit branched covers of degree at most $60$
which in turn are also branched covers of the standard spherical football.

One may wonder what special r\^ole pentagons and hexagons play in this discussion. We shall address this question 
in the second part of this paper, where we determine all the possibilities for triples $(k,l,n)$ that can be realised by 
maps on the two-sphere whose faces are $k$-gons and $l$-gons satisfying the conditions that the edges of each $k$-gon 
meet only edges of $l$-gons, and that every $n^{\textrm{th}}$ edge of each $l$-gon meets an edge of a $k$-gon, and its 
other edges meet $l$-gons. Not surprisingly, the determination of these triples is a generalisation of the topological argument 
determining the Platonic solids. In most cases the generalised football patterns have realisations dual to very symmetric 
triangulations of the sphere that have been known and studied since the $19^{\textrm{th}}$ century. That purely 
topological or combinatorial considerations lead to a list that contains almost only the usual symmetric patterns and 
their degenerations is a kind of rigidity phenomenon associated with these spherical triangulations.

We shall see that the classification theorem for spherical football patterns proved in the first part of this paper for the triple 
$(5,6,2)$ actually holds for all generalised football patterns with $n=2$: each generalized football with a pattern of type 
$(k,l,2)$ is a branched cover of the corresponding minimal pattern. We shall also see that this result does not extend to $n>2$.


\subsubsection*{Acknowledgement}

We are grateful to B.~Hanke and B.~Leeb for helpful comments, and to the {\it Deutsche Forschungsgemeinschaft} for support.

\section{Ribbon graphs and branched covers}

\subsection{Football graphs}

A football pattern is a map in the sense of~\cite{Maps} on the two-sphere satisfying the usual conditions that at least three 
edges meet at every vertex, that all faces are pentagons and hexagons, that the edges of each pentagon meet only edges 
of hexagons, and that the edges of each hexagon alternately meet edges of pentagons and of hexagons. We make no 
regularity assumption, so that a football pattern is not a geometric, but a combinatorial-topological object.

A football pattern determines, and is determined by, its dual graph. This graph has a vertex for every polygon in the pattern, 
and the vertices are coloured, say black for the vertices corresponding to pentagons and white for the vertices corresponding 
to hexagons. Two vertices are connected by an edge if the corresponding polygons share an edge. The edges meeting at a 
vertex are cyclically ordered (with respect to this endpoint) by remembering that the sides of a polygon are cyclically ordered. 
Therefore the dual graph is a fatgraph or ribbon graph in the sense of~\cite{RibbonGraphs}, leading to the following definition:
\begin{defn}\label{d:basic}
A football graph is a ribbon graph with black and white vertices satisfying the following conditions:
\begin{enumerate}
\item each black vertex has valence five, and all five edges connect the given vertex to white vertices, and
\item each white vertex has valence six, and the six edges alternately\footnote{The alternating condition is with respect to 
the cyclic order of the edges.} connect the given vertex to black and white vertices.
\end{enumerate}
\end{defn}
The picture at the beginning of this paper shows the dual ribbon graph of the standard football pattern on $S^2$ consisting of $12$
pentagons and $20$ hexagons. To make the picture symmetric, a black vertex is spread out at infinity and is depicted by
the shaded ring around the rest of the graph.

As every finite ribbon graph corresponds to a unique closed oriented surface, we have a natural bijection between football
graphs and football patterns on arbitrary closed oriented surfaces. A covering map between football graphs corresponds to 
a possibly branched covering map between surfaces, with any branching restricted to the centers of the faces of 
the decompositions given by the football graphs. As the football graphs are dual to actual football patterns, branching
can only occur at the vertices of patterns.

Let $b$ and $w$ be the numbers of black and white vertices in a football graph $\Gamma$. 
Then the total number of vertices is $v=b+w$, and the number of edges is $e=\frac{1}{2}(5b+6w)$.
Of these edges, $e_1=\frac{3}{2}w$ have white endpoints, and $e_2=5b=3w$ have a black and a white 
endpoint. It follows that there is a natural number $d$ such that $b=6d$ and $w=10d$.

\begin{lem}\label{l:triangle}
The only football graph giving a triangulation of a closed surface is the dual graph $\Gamma_0$ of the standard spherical football
with $d=2$.
\end{lem}
\begin{proof}
Let $\Sigma$ be the closed oriented surface defined by a football graph $\Gamma$. By Euler's formula, the number of faces in
the cell decomposition of $\Sigma$ determined by $\Gamma$ is
$$
f = \chi (\Sigma ) + e - v =  \chi (\Sigma ) + 45d - 16d =  \chi (\Sigma ) + 29d  \ .
$$
If $\Gamma$ defines a triangulation of $\Sigma$, then we must have $2e = 3f$, which, rewritten in terms of $d$, means
$90d = 3\chi (\Sigma )+ 87d$, or $d=\chi (\Sigma )$. Thus $d=2$, and $\Sigma$ is $S^2$. The combinatorics of the 
corresponding football pattern is uniquely determined in this case, as can be seen from the  proof of Theorem~\ref{t:main}
below.
\end{proof}

\subsection{The football tree}

There is precisely one connected and simply connected football graph, which we shall call the football tree $T$. It is the universal 
cover of any football graph. As $T$ thought of as a ribbon graph is an orientable surface, it makes sense to speak of 
orientation-preserving automorphisms, and we shall 
denote the group of all such automorphisms by $\Aut(T)$. This can be determined explicitly:
\begin{thm}
The automorphism group $\Aut(T)$ of the football tree is isomorphic to the free product $\Z_2\star\Z_3\star\Z_5$.
\end{thm}
\begin{proof}
This is a straightforward application of the Bass--Serre theory of groups acting on trees, cf.~\cite{Trees}. This theory is usually
formulated for groups acting on trees without inverting edges. In our situation, there are automorphisms inverting edges 
that connect a pair of white vertices. Therefore, we subdivide each of these edges by introducing red vertices in the 
middle of each edge connecting two white edges of $T$. We obtain a new tree $T'$, which has three kinds of edges: the white 
and black ones of valence $6$ and $5$ respectively, and red ones of valence $2$. The red and black vertices are only connected 
to white ones, and the edges meeting at a white vertex lead alternately to red and black vertices. Now $\Aut(T)$ acts on $T'$
without inverting any edges. The action is simply transitive on the two kinds of edges, black-white and red-white. The action is 
also transitive on the vertices of a given colour, with isotropy groups of orders $2$, $3$ and $5$ for the red, white and black 
vertices. The quotient graph $T'/\Aut(T)$ is a tree with three vertices, one for each colour, and with one edge connecting the 
white vertex to each of the other vertices. We think of this as a graph of groups by labeling the vertices with the isotropy groups. 
As the edges have trivial isotropy, the fundamental group of this graph of groups is the free product of the labels of the vertices.
By the structure theorem of~Section~I.5.4 in~\cite{Trees}, this fundamental group is isomorphic to $\Aut(T)$.
\end{proof}

\subsection{The classification of spherical footballs}

We now want to prove that every football pattern on the sphere is obtained as a branched cover of the standard football
branched only in the vertices. Equivalently, we prove that the spherical dual ribbon graphs are all obtained as covering spaces of 
the dual graph of the standard football. The first step is the following:
\begin{lem}\label{l:prep}
Let $\Gamma$ be a football graph with universal covering $\pi\colon T\longrightarrow\Gamma$. Suppose that $\gamma$ is an
oriented path in $T$ consisting of a sequence of edges without backtracking. If $\pi$ maps $\gamma$ to a closed path that corresponds
to a boundary component of $\Gamma$ thought of as a surface with boundary, then $\gamma$ consists of $3n$ edges for some 
natural number $n$, and the image of $\gamma$ in the standard football graph $\Gamma_0$ is a loop that is the $n^{\textrm{th}}$
power of the loop formed by a triangle in the triangulation defined by $\Gamma_0$, cf.~Lemma~\ref{l:triangle}.
\end{lem}
\begin{proof}
We think of $T$ as a surface with boundary.
Choose a boundary component $C$ covering the boundary component $\pi(C)$ of $\Gamma$ to which $\gamma$ is mapped.
As a boundary component of $T$, $C$ runs right along a sequence $\{e_i\}_{i\in\ZZ}$ of oriented edges in $T$ such that, of course, 
the origin $o(e_i)$ of each edge coincides with the endpoint of the previous edge, and, in addition, with respect to the cyclic order of 
edges emanating from $o(e_i)=o(\bar e_{i-1})$, $e_i$ is the successor of $\bar e_{i-1}$. (The bar denotes edge inversion.) It follows 
from Definition~\ref{d:basic} that the sequence of vertices $o(e_i)$ is of the form \emph{black, white, white, black, white, white}, etc.

The setwise stabilizer of $C$ in $\Aut(T)$ is the infinite cyclic group generated by
the translation $\tau$, which maps $e_i$ to $e_{i+3}$.  Now $\pi(C)=C/G$, for some subgroup $G \le \Stab(C)$,
i.\,e.\ $G = \langle \tau^n \rangle$ for some $n\in\NN$.  Thus $\pi(C)$ runs
along a $3n$-gon, which must be $\pi(\gamma)$.

By assumption, the path $\gamma$ is a piece $e_{j+1},\ldots,e_{j+k}$ of the sequence $\{e_i\}$ and projects to 
$\pi(\gamma)$. Therefore $k=3n$. The images of $e_i$ in the standard football still satisfy the condition that 
consecutive edges be related by the cyclic ordering defining the ribbon graph structure.
But this means that the image of $e_{j+1}e_{j+2}e_{j+3}$ is a triangle.
\end{proof}
Here is the classification theorem for spherical football patterns:
\begin{thm}\label{t:main}
Every football graph dual to a football pattern on $S^2$ is a finite covering space of the standard football graph $\Gamma_0$.
Equivalently every football pattern on $S^2$ is obtained from the standard one by passing to a branched cover branched only in
vertices of the pattern.
\end{thm}
\begin{proof}
Given a spherical football graph $\Gamma$, fix a universal covering $\pi\colon T\rightarrow\Gamma$ of $\Gamma$ and a 
universal covering $\pi_0\colon T \rightarrow \Gamma_0$ of the standard football graph $\Gamma_0$.  We are going to show 
that the group of deck transformations $\Aut_\Gamma T$ is a subgroup of $\Aut_{\Gamma_0}T$.  This implies that $\pi_0$ 
factors through $\pi$.

Choose a point $p$ on an edge of $T$ that is not an endpoint or a midpoint, so that it has trivial stabilizer in $\Aut (T)$.  
Covering space theory identifies the group of deck transformations $\Aut_{\Gamma}T$ with the fundamental group 
$\pi_1(\Gamma;\pi(p))$.  Since $\Gamma$ is spherical, this fundamental group is generated by paths of the form 
$\beta\gamma\beta^{-1}$, where $\gamma$ is a loop along a boundary component, $\beta$ runs from the base point $\pi(p)$
to the origin of $\gamma$ and $\beta^{-1}$ is the way back.  Lifting $\beta$ and $\gamma$ to $T$ we obtain a path
$\tilde\beta\tilde\gamma\tau^n(\tilde\beta^{-1})$ leading from $p$ to $\tau^n(p)$, with $\tau$ from the proof of the previous lemma.  
Hence $\tau^n$ is the deck transformation corresponding to the given generator of $\pi_1(\Gamma;\pi(p))$.  This proves the result, 
because $\tau$ is a deck transformation over $\Gamma_0$.
\end{proof}
This proof also shows:
\begin{cor}
The subgroup $\pi_1(\Gamma_0)\subset\Aut(T)$ is normal. The quotient $\Aut(T)/\pi_1(\Gamma_0)$ is the icosahedral group
of order $60$.
\end{cor}

\subsection{Footballs of positive genera}

Every ribbon graph corresponds to a unique closed oriented surface, and of course every
such surface does indeed arise from a football graph, for example because it is a branched covering
of the two-sphere, which we  can arrange to be branched only in the vertices of a suitable
football pattern. We now want to show that there are other football patterns on surfaces of 
positive genera, that are not lifted from the two-sphere. The proof of Theorem~\ref{t:main}
does not extend, because for a ribbon graph corresponding to a surface of positive genus 
there are generators in the fundamental group that arise from handles, rather than from the 
punctures.

Recall that the parameter $d$ for a finite football graph specifies the number of black and white
vertices by the formulae $b=6d$ and $w=10d$. Passing to a $D$-fold covering multiplies $d$ by
$D$. As the standard football graph $\Gamma_0$ has $d=2$, all its non-trivial coverings have $d\geq 4$.
Therefore, to exhibit football graphs that are not coverings of $\Gamma_0$, it suffices to find
examples of positive genus with $d<4$. Performing certain cut-and-paste operations on $\Gamma_0$,
we can actually produce examples with $d=2$ and rather large genera.

The simplest example is the following.
\begin{ex}\label{ex:sum}
Pick two disjoint edges in the standard football pattern, that are of the same type, so that they both 
separate hexagons from each other, or they both separate a pentagon from a hexagon. Open up
the two-sphere along these edges to obtain a cylinder whose two boundary circles each have two vertices 
and two edges. As the  two edges along which we opened the sphere were of the  same type, the 
two boundary circles of the cylinder can be identified in such a way that the resulting torus carries
an induced football pattern with $d=2$. 
\end{ex}
In this example there are $58$ vertices instead of the $60$ in the standard spherical football. All but
two of them are $3$-valent, and the exceptional two are $6$-valent.

In the language of  ribbon graphs, the surgery performed in the above example amounts to cutting two
ribbons and regluing the resulting ends in a different pairing. This can also be done with ribbons 
corresponding to edges that share a vertex, in which case instead of cutting and pasting, the surgery
can be described through the reordering of edges:
\begin{ex}\label{ex:surgery}
Let $e_1,\ldots,e_5$ be the edges emanating from a black vertex in the standard football graph $\Gamma_0$.  
Define a new football graph $\Gamma$ by reordering the edges as $e_1$, $e_3$, $e_2$, $e_4$, $e_5$.  
This procedure glues the three triangles whose edges include $e_2$ or $e_3$ into a single $9$-gon 
boundary component of $\Gamma$.  This new graph still has $d=2$, but the underlying surface is a torus.
In the dual football pattern there are $58$ vertices, of which $57$ are $3$-valent and one is $9$-valent.
\end{ex}

\begin{ex}\label{ex:symmetric-genus-24}
Let $e_1,\ldots,e_6$ be the edges emanating from a white vertex in the standard football graph $\Gamma_0$, 
enumerated in their cyclic order and labelled such that $e_1$, $e_3$ and $e_5$ have black ends.  Define a new 
ribbon graph by reordering the edges cyclically as $e_1,e_4,e_3,e_6,e_5,e_2$. This means that the edges 
leading to white vertices are cut and reattached after a cyclic permutation given geometrically by a rotation 
by angle $\frac{2\pi}{3}$ around the vertex in the realisation of $\Gamma_0$ with icosahedral symmetry.

We apply this procedure to every white vertex of $\Gamma_0$. The resulting football graph $\Gamma$ is 
symmetric in the sense that it admits rotations of order $2$, $3$ and $5$ around an edge, a white vertex and 
a black vertex, respectively.  Hence the full football group $\Aut(T)$ acts on $\Gamma$.  In particular all faces 
are conjugate.  One can verify by inspection that the faces are $15$-gons.  Hence the Euler characteristic of the 
underlying surface is $-46$, and its genus is $24$. This is a football graph with $d=2$ and is therefore not a 
covering of $\Gamma_0$.
\end{ex}

Although football graphs of positive genera are not in general coverings of $\Gamma_0$, we have the 
following:
\begin{thm}
Every football pattern on a closed oriented surface admits a branched cover of degree $D\leq 60$ that is also
a branched cover of the standard minimal pattern on $S^2$. The bound for $D$ is sharp.
\end{thm}
\begin{proof}
Let $\Gamma$ be a finite football graph. As $\pi_1(\Gamma_0)\subset\Aut(T)$ is a subgroup of index $60$,
the intersection $\pi_1(\Gamma_0)\cap\pi_1(\Gamma)$ has index at most $60$ in $\pi_1(\Gamma)$. The 
intersection corresponds to a covering of $\Gamma$ of degree $D\leq 60$ that is also a covering of $\Gamma_0$.

To prove that coverings of degree strictly less than $60$ do not always suffice, recall that $\Aut(T)$ acts on the 
genus $24$ football graph $\Gamma$ from Example~\ref{ex:symmetric-genus-24}. This is equivalent to the
fundamental group $\pi_1(\Gamma)$ being a normal subgroup of $\Aut(T)$. The covering of $\Gamma$ corresponding 
to the subgroup $N=\pi_1(\Gamma_0)\cap\pi_1(\Gamma)$ of $\Aut(T)$ is a Galois covering with Galois group 
$\pi_1(\Gamma)/N$. Now the injection $\pi_1(\Gamma)\map\Aut(T)$ induces
an embedding of $\pi_1(\Gamma)/N$ as a normal subgroup of the
orientation-preserving icosahedral group $\Aut(T)/\pi_1(\Gamma_0)$,
isomorphic to the alternating group $A_5$.  Since this is a simple group we
must have $\pi_1(\Gamma)/N \cong A_5$ or $\{1\}$.  The second case is excluded
because $\Gamma$ is not a covering of $\Gamma_0$.
\end{proof}

\subsection{Non-orientable footballs}

Although we have modelled football patterns as ribbon graphs, we can also consider them on non-orientable 
surfaces, because the condition that every other edge emanating from a white vertex should connect to a 
black vertex is preserved by inversion of the cyclic order.

Here are the simplest examples for the projective plane.
\begin{ex}
A football pattern on the real projective plane is readily constructed from the standard football.  In the 
dual ribbon graph $\Gamma_0$ cut a single ribbon and reglue it with a half-twist so that the surface becomes 
non-orientable. This gives a football pattern with $d=\text2$ that, instead of the $60$ vertices of valence $3$
in the standard football, has $58$ vertices of valence $3$ and a unique vertex of valence $6$. Therefore the 
Euler number of the surface is $1$.
\end{ex}
If we lift the pattern in this example to the universal covering of the projective plane, we obtain a football pattern
with $d=4$ on the two-sphere, which, by Theorem~\ref{t:main}, is a $2$-fold branched cover of the standard pattern.
Of course in this case the branched covering structure can be seen directly, by focussing on the two vertices of 
valence $6$.

\begin{ex}
The usual symmetric realisation of the standard football pattern on the sphere as a truncated icosahedron is 
symmetric under the antipodal involution.  Thus it descends to a pattern on $\R P^2$ with $d=\text1$.
\end{ex}

\section{Generalised football patterns}

In this section we consider generalisations of the traditional football patterns.
\begin{defn}
A generalised football pattern is a map on the two-sphere whose faces are $k$-gons and $l$-gons satisfying the conditions 
that the edges of each $k$-gon meet only edges of $l$-gons, and that every $n^{\textrm{th}}$ edge of each $l$-gon meets 
an edge of a $k$-gon, and its other edges meet $l$-gons. 
\end{defn}
To avoid degenerate cases we always assume $k\geq 3$, $l\geq 3$ and $l=n\cdot m$ with positive integers $m$ 
and $n$. As usual, at least three edges meet at every vertex.

If a given triple $(k,l,n)$ can be realised by a generalised football pattern, then it can be realised in infinitely many ways, for
example by taking branched covers branched only in the vertices of a given pattern. We will determine all possible triples,
and we will find a minimal realisation for each of them. We will also see that in some cases there are realisations that
are not branched covers of the minimal one.

Before proceeding to the classification, we list some examples for future reference.

\subsection{Some examples}

The standard football realising the triple $(5,6,2)$ can be thought of as a truncated icosahedron. More generally, we have:
\begin{ex}\label{ex:tP}
The truncated Platonic solids realise the triples $(3,6,2)$, $(3,8,2)$, $(3,10,2)$, $(4,6,2)$ and $(5,6,2)$ as generalised 
football patterns.
\end{ex}
There are also infinite series of examples obtained by truncating the degenerate Platonic solids:
\begin{ex}\label{ex:Amfoot}
Start with a subdivision of the sphere along $k\geq 3$ halves of great circles running from the north to the south poles.
We shall call this an American football. If we now truncate at one of the poles, we obtain a realisation of the triple
$(k,3,3)$. (For $k=3$ this is a tetrahedron.) If we truncate at both poles, we obtain a tin can pattern realising $(k,4,2)$. 
(For $k=4$ it is a cube.) This is also known as a $k$-prism. If we add $k$ edges along the equator to
this last example, we obtain a double tin can realising $(k,4,4)$, for any $k\geq 3$.
\end{ex}
Here is a variation on the above tin can pattern:
\begin{ex}\label{ex:zz}
Take a $k$-gon, with $k\geq 3$ arbitrary, and surround it by pentagons in such a way that the two pentagons meeting 
a pair of adjacent sides of the $k$-gon share a side. We can fit together two such rings made up of a $k$-gon and 
$k$ pentagons each along a zigzag curve to obtain a realisation of the triple $(k,5,5)$ by a generalised football pattern.
(For $k=5$ we obtain a dodecahedron.)
\end{ex}
The next two examples need to be visualised using the accompanying figures.
\begin{ex}\label{ex:3}
Take a Platonic solid, and subdivide each face as follows. If the face is a $k$-gon, put a smaller $k$-gon in its interior, and 
radially connect each corner of this smaller $k$-gon to the corresponding corner of the original face. In this way each face
of the original polyhedron is divided into a $k$-gon and $k$ quadrilaterals. The cases $k=3$ and $4$ are shown in Figure~\ref{figure1}. 
\begin{figure}
\epsfig{file=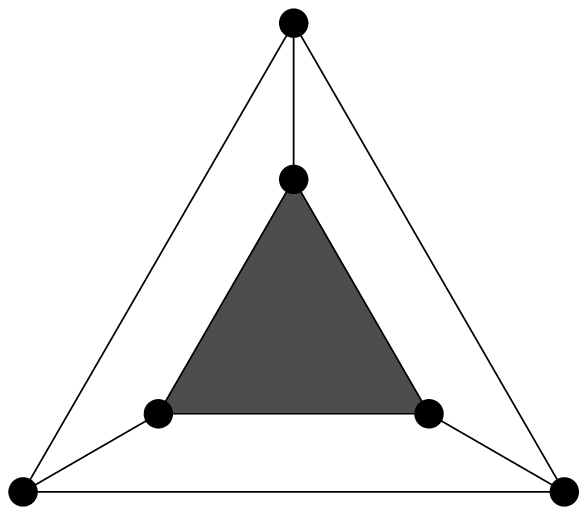}
\epsfig{file=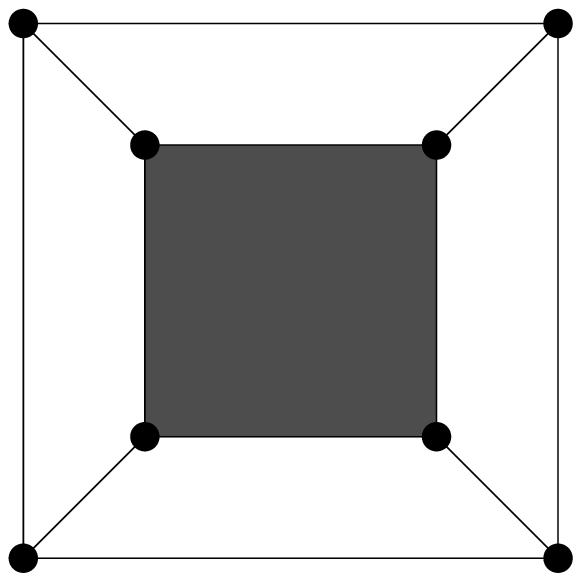}
\caption{Variation on Platonic solids}\label{figure1}
\end{figure}
Next we erase the edges of the original polyhedron, so that the two quadrilaterals of the subdivison meeting along an edge 
are joined to form a hexagon. In this way we obtain a realisation of the triple $(k,6,3)$. We shall refer to this construction as 
a variation on the original Platonic solid.
\end{ex} 
\begin{ex}\label{ex:6}
Again we start with a Platonic solid whose faces are $k$-gons. We subdivide each face into a smaller $k$-gon and $k$
hexagons as shown in Figure~\ref{figure2}. 
\begin{figure}
\epsfig{file=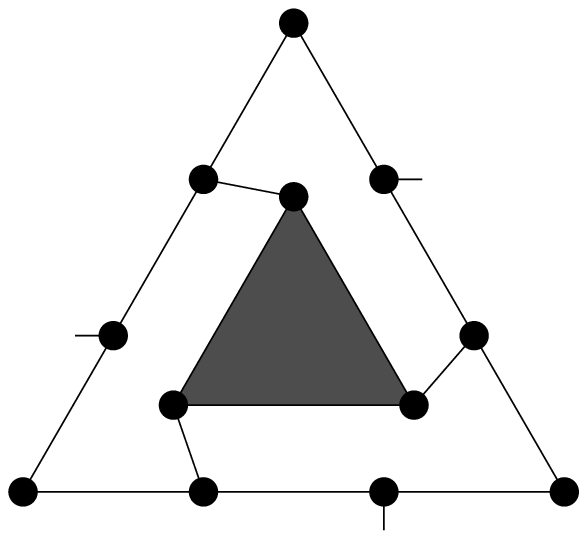}
\epsfig{file=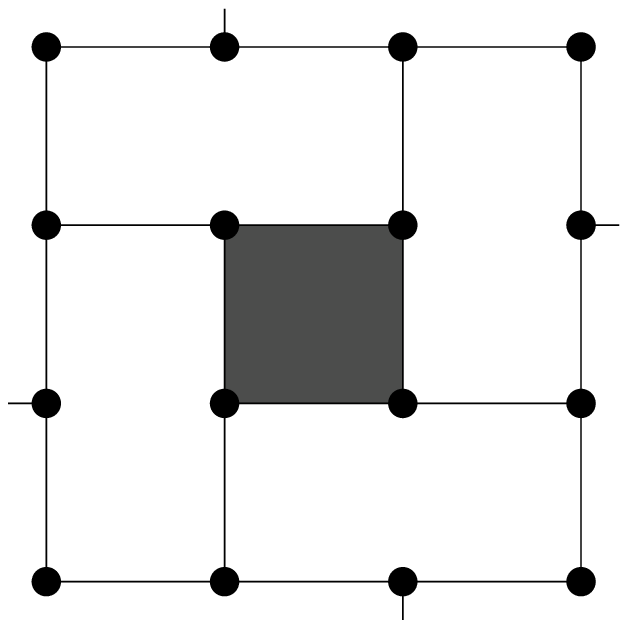}
\caption{Subdivision of Platonic solids}\label{figure2}
\end{figure}
This gives a realisation of the triple $(k,6,6)$.
\end{ex}

\subsection{The classification of generalised football patterns}

Now we want to prove that the previous examples exhaust all possible generalised football patterns with $n\geq 2$. We shall 
also treat the case $n=1$. The results of this classification are summarised in the table in Figure~\ref{fig:table}.

If $S^2$ is endowed with a generalised football pattern of type $(k,l,n)$, we shall think of the $k$-gons as being black and the 
$l$-gons as being white. Their numbers are denoted by $b$ and $w$ respectively. The pattern then has $f=b+w$ many faces, 
it has $e=\frac{1}{2}(bk+wl)$ many edges, and the number $v$ of vertices is bounded by $v\leq\frac{1}{3}(bk+wl)$ as there 
have to be at least $3$ faces meeting at every vertex.

Counting the number of edges at which a $k$-gon meets an $l$-gon in two different ways, we find
\begin{equation}\label{eq:fund}
k\cdot b=\frac{1}{n}\cdot l\cdot w = m\cdot w \ .
\end{equation}
Computing the Euler characteristic and using~\eqref{eq:fund} leads to:
\begin{equation*}
2=f-e+v\leq b+w-\frac{1}{6}(k\cdot b+l\cdot w)=b+w-\frac{1}{6}k\cdot b\cdot (n+1) \ .
\end{equation*}
Dividing by $2k\cdot b$ we obtain the key inequality
\begin{equation}\label{eq:key}
\frac{1}{k\cdot b}+\frac{n+1}{12}\leq \frac{1}{2k}+\frac{1}{2m} \ .
\end{equation}

\subsubsection{The case $n\geq 2$}

If $n\geq 2$, then the left hand side of~\eqref{eq:key} is strictly larger than $\frac{1}{4}$. If both $k$ and $m$ are at least $4$, 
then the right hand side is at most $\frac{1}{4}$. Thus $n\geq 2$ implies $k=3$ or $m\in\{1,2,3\}$. We now discuss these cases 
separately. 

For every possible triple $(k,l,n)$ the inequality~\eqref{eq:key} gives a lower bound for $b$. With~\eqref{eq:fund} we then have 
a lower bound for $w$. Checking these values against the examples of the previous section, one sees for most of the examples 
that they actually give minimal realisations. There are only a few cases when this simple check does not suffice, because the minimal
realisations have vertices of valence strictly larger than $3$.


\begin{lem}\label{l:k=3}
Suppose the triple $(k,l,n)$ with $k=3$ and $n\geq 2$ is realised by a generalised football pattern. Then $(k,l,n)$ is one of the 
triples $(3,4,2)$, $(3,6,2)$, $(3,8,2)$, $(3,10,2)$, $(3,3,3)$, $(3,6,3)$, $(3,4,4)$, $(3,5,5)$, $(3,6,6)$. Minimal realisations are 
given by the Examples~\ref{ex:tP}, 
\ref{ex:Amfoot}, \ref{ex:zz}, \ref{ex:3} and~\ref{ex:6}.
\end{lem}
\begin{proof}
Putting $k=3$ in~\eqref{eq:key} we obtain
\begin{equation}\label{eq:3key}
\frac{1}{3b}+\frac{n}{12}\leq \frac{1}{12}+\frac{1}{2m} \ .
\end{equation}
Together with $n\geq 2$, this implies $m\leq 5$.

If $m=5$, then we obtain $n=2$ from~\eqref{eq:3key}. Thus $l=10$.
The truncated dodecahedron from Example~\ref{ex:tP} is a minimal realisation.

If $m=4$, then again $n=2$ from~\eqref{eq:3key}. Thus $l=8$.
The truncated cube from Example~\ref{ex:tP} is a minimal realisation.

If $m=3$, then again $n=2$ from~\eqref{eq:3key}. Thus $l=6$.
The truncated tetrahedron from Example~\ref{ex:tP} is a minimal realisation.

If $m=2$, then~\eqref{eq:3key} gives $n\leq 3$. 
For $n=3$, equivalently $l=6$, a minimal realisation is the variation on the tetrahedron from Example~\ref{ex:3}. 
For $n=2$, equivalently $l=4$, we get the case $k=3$ of the truncated American football in Example~\ref{ex:Amfoot}.

If $m=1$, then~\eqref{eq:3key} only gives $n\leq 6$.
For $n=6$, therefore $l=6$, the subdivision of the tetrahedron from Example~\ref{ex:6} gives a minimal realisation.
For $n=l=5$, the case $k=3$ in Example~\ref{ex:zz} is a minimal realisation.
For $n=l=4$, the case $k=3$ of the double tin can in Example~\ref{ex:Amfoot} is a realisation.
As it has vertices of valence $4$, it is not immediately obvious that it is a minimal realisation. In this case~\eqref{eq:3key}
gives $b\geq 1$, but one can easily check that $b=1$ is not possible. Thus $b$ is at least $2$, and the double tin can
is indeed a minimal realisation.
For $n=l=3$, the case $k=3$ of the partially truncated American football in Example~\ref{ex:Amfoot} is a minimal realisation.
Note that $m=1$ implies $n\geq 3$, so that we have now exhausted all cases with $k=3$ and $n\geq 2$.
\end{proof}


\begin{lem}\label{l:m=3}
Suppose the triple $(k,l,n)$ with $l/n=m=3$ and $n\geq 2$ is realised by a generalised football pattern. Then $(k,l,n)$ is one of the 
triples $(3,6,2)$, $(4,6,2)$, $(5,6,2)$.
Minimal realisations are given by the truncated Platonic solids in Example~\ref{ex:tP}.
\end{lem}
%
%
%


\begin{lem}\label{l:m=2}
Suppose the triple $(k,l,n)$ with $l/n=m=2$ and $n\geq 2$ is realised by a generalised football pattern. Then $(k,l,n)$ is one of the 
triples $(3,6,3)$, $(4,6,3)$, $(5,6,3)$, or $(k,4,2)$ with $k\geq 3$.
Minimal realisations are given by the variations on the Platonic solids in Example~\ref{ex:3}, respectively by the truncated American
football in Example~\ref{ex:Amfoot}.
\end{lem}
%
%
We omit the proofs of Lemmas~\ref{l:m=3} and~\ref{l:m=2}, because they are completely analogous to, and even simpler than, the 
proof of Lemma~\ref{l:k=3}.


\begin{lem}\label{l:m=1}
Suppose the triple $(k,l,n)$ with $l/n=m=1$ and $n\geq 2$ is realised by a generalised football pattern. Then $(k,l,n)$ is one of the 
triples $(3,6,6)$, $(4,6,6)$, $(5,6,6)$, or $(k,3,3)$, $(k,4,4)$, $(k,5,5)$ with $k\geq 3$.
Minimal realisations are given by the subdivisions on the Platonic solids in Example~\ref{ex:6}, respectively by the infinite sequences
in Examples~\ref{ex:Amfoot} and~\ref{ex:zz}.
\end{lem}
\begin{proof}
For $m=1$, equivalently $l=n$, we obtain $n\leq 6$ for all $k\geq 3$ from~\eqref{eq:key}.

If $l=n=6$, then we also have $k\leq 5$ from~\eqref{eq:key}. Thus $k$ is $3$, $4$ or $5$, and minimal realisations are given in 
Example~\ref{ex:6}.

If $l=n=5$, then all $k\geq 3$ are possible, and minimal realisations are given in Example~\ref{ex:zz}.

If $l=n=4$, then all $k\geq 3$ are possible, and realisations are given by the double tin cans in Example~\ref{ex:Amfoot}.
To see that these realisations are minimal, it suffices to exclude the case $b=1$, which is easily done by contradiction.

If $l=n=3$, then again all $k\geq 3$ are possible, and minimal realisations are given by the partially truncated American football in 
Example~\ref{ex:Amfoot}.
\end{proof}

This completes the classification of generalised football patters with $n\geq 2$ on the two-sphere.

\subsubsection{The case $n=1$}

A generalised football pattern with $n=1$ consists of $b$ black $k$-gons and $w$ white $l$-gons so that the two polygons
meeting along an edge always have different colours. Note that here the situation is completely symmetric in $k$ and $l$.

\begin{lem}\label{l:n=1}
Suppose the triple $(k,l,1)$ is realised by a generalised football pattern on $S^2$. Then, up to changing the r\^oles of $k$ and $l$,
$(k,l)$ is one of the pairs $(3,3)$, $(3,4)$ or $(3,5)$. Minimal realisations are obtained by painting the faces of an octahedron,
a cuboctahedron, respectively an icosidodecahedron, in a suitable manner.
\end{lem}
\begin{proof}
Counting the edges in two different ways leads to $b\cdot k = w\cdot l$. 
As every edge separates a black from a white polygon, there must be an even number of faces meeting at every vertex. 
Therefore the valence of every vertex is $\geq 4$, giving rise to
$$
v\leq\frac{1}{4}(b\cdot k + w\cdot l) \ ,
$$
which is of course stronger than what we had before, when the valence of a vertex was only $\geq 3$.
Computing the Euler characteristic as before, we  obtain instead of~\eqref{eq:key} the stronger
\begin{equation}\label{eq:key1}
\frac{1}{k\cdot b}+\frac{1}{4}\leq \frac{1}{2k}+\frac{1}{2l} \ .
\end{equation}
If both $k$ and $l$ are $\geq 4$, then the right hand side is $\leq\frac{1}{4}$, which is impossible. Thus $k$ or $l$ is
$=3$. By the symmetry between $k$ and $l$ we may assume that $k=3$. Substituting this into~\eqref{eq:key1}, we 
find
\begin{equation}\label{eq:key2}
\frac{1}{3b}+\frac{1}{12}\leq \frac{1}{2l} \ ,
\end{equation}
which implies $l\leq 5$.

If $l=3$, then $b=w\geq 4$. A realisation is obtained by painting the faces of an octahedron in black and white, so that
each edge separates a black face from a white one.

If $l=4$, then~\eqref{eq:key2} implies $b\geq 8$, so that $w\geq 6$. The cuboctahedron, cf.~\cite{S}, is a realisation.

If $l=5$, then~\eqref{eq:key2} implies $b\geq 20$, so that $w\geq 12$. The icosidodecahedron, cf.~\cite{S}, is a realisation.

All these realisations are minimal, because they have precisely four faces meeting at every vertex, so that~\eqref{eq:key2}
becomes an equality.
\end{proof}

We summarise the above classification of the generalised football patterns as follows:
\begin{thm}\label{t:table}
Suppose that $S^2$ admits a map whose faces are $k$-gons and $l$-gons with $k,l\geq 3$, such that the edges of
each $k$-gon meet only edges of $l$-gons, and so that every $n^{\textrm{th}}$ edge of each $l$-gon meets an edge of a $k$-gon, 
and its other edges meet $l$-gons. Then $l\leq 10$ and $n\leq 6$. There are $16$ different sporadic triples $(k,l,n)$ with $k\leq 5$ that 
occur, together with $4$ infinite sequences with variable $k$ and fixed $l$ and $n$. All the possibilities are listed in the table in Figure~\ref{fig:table},
which also gives minimal realisations for all cases.
\end{thm}
\begin{figure}
$$ 
\vbox{\tabskip=.5em\offinterlineskip
\halign{\strut#&\hfil#\hfil & \vrule# & \hfil#\hfil & \vrule# &
\hfil#\hfil & \vrule# & \hfil#\hfil & \vrule# & \hfil#\hfil & \vrule# &
\hfil#\hfil & \vrule# & \hfil#\hfil\cr
&  && k && m && n && {\bf minimal realisation} && b && w \cr
\noalign{\hrule} 
\noalign{\hrule} 
& 1. && 3 && 3 && 1 && octahedron && 4 && 4 \cr
& 2. && 3 && 4 && 1 && cuboctahedron && 8 && 6 \cr
& 3. && 4 && 3 && 1 && cuboctahedron && 6 && 8 \cr
& 4. && 3 && 5 && 1 && icosidodecahedron && 20 && 12 \cr
& 5. && 5 && 3 && 1 && icosidodecahedron && 12 && 20 \cr
& 6. && 3 && 3 && 2 && truncated tetrahedron && 4 && 4 \cr
& 7. && 3 && 4 && 2 && truncated cube && 8 && 6 \cr
& 8. && 4 && 3 && 2 && truncated octahedron && 6 && 8 \cr
& 9. && 3 && 5 && 2 && truncated dodecahedron && 20 && 12 \cr
& 10. && 5 && 3 && 2 && truncated icosahedron = {\bf football} && 12 && 20 \cr
& 11. && $\geq$ 3 && 2 && 2 && truncated American football && 2 && k \cr
& 12. && 3 && 2 && 3 && variation on the tetrahedron && 4 && 6 \cr
& 13. && 4 && 2 && 3 && variation on the cube && 6 && 12 \cr
& 14. && 5 && 2 && 3 && variation on the dodecahedron && 
12 && 30 \cr
& 15. && $\geq$ 3 && 1 && 3 && partially truncated American football && 1 && 
k \cr
& 16. && $\geq$ 3 && 1 && 4 && double tin can && 2 && 2k \cr
& 17. && $\geq$ 3 && 1 && 5 && zigzag tin can && 2 && 2k \cr
& 18. && 3 && 1 && 6 && subdivision of the tetrahedron && 4 && 12 \cr
& 19. && 4 && 1 && 6 && subdivision of the cube && 6 && 24 \cr
& 20. && 5 && 1 && 6 && subdivision of the dodecahedron && 12 && 60 \cr
}} 
$$
\caption{The classification of generalised football patterns on $S^2$}\label{fig:table}
\end{figure}
The different minimal realisations were described in our earlier examples, and in the course of the proof.
There are several alternative descriptions of items 12.--14. The variation on the tetrahedron is nothing but a partially truncated
cube, truncated at $4$ of its $8$ vertices, chosen so that each face is truncated at two diagonally opposite corners. Similarly,
the variations on the cube and the dodecahedron are partial truncations of the rhombic dodecahedron and of the rhombic 
triacontrahedron respectively.
The subdivisions in items 18.--20.~can also be thought of as partial truncations of the 
dodecahedron, the pentagonal icositetrahedron, respectively the pentagonal hexecontrahedron.

In order to make the symmetries more obvious, the table does not list triples $(k,l,n)$, but rather
$(k,m,n)$ with $m=l/n$. Note that this makes no difference when $n=1$, in which case there is a complete
symmetry between $k$ and $l=m$. Therefore, the cases 2.~and 3., respectively 4.~and 5., are 
dual to each other, obtained by switching the roles of $k$ and $l$. Case 1.~is self-dual.
Similarly, cases 7.~and 8., respectively 9.~and 10., are dual to each other with the duality 
induced by the duality of Platonic solids, and case 6.~is self-dual.

\subsection{Branched covers for generalised football patterns}

Now that we have an overview of all the generalised football patterns, one may ask whether Theorem~\ref{t:main}
can be extended, to prove that every generalised spherical football is a branched cover of the corresponding minimal 
example, branched only in the vertices. It is not hard to see that for $n=2$ this is indeed the case:

\begin{thm}\label{t:cover}
If $(k,l,2)$ is realised by a generalised football pattern on $S^2$, then every spherical realisation is a branched 
cover, branched only in the vertices, of the minimal realisation. In particular the minimal realisation is unique.
\end{thm}
\begin{proof}
The minimal realisations for $n=2$ have all vertices of valence $3$. Equivalently, the dual graph triangulates the 
sphere. Moreover, every realisation has the property that the valence of every vertex is a multiple of $3$, with
every third face a black $k$-gon. Therefore the proofs of Lemma~\ref{l:prep} and of Theorem~\ref{t:main} go through.
\end{proof}
For larger values of $n$ one loses control of the structure of the vertices, and the proof breaks down.
\begin{ex}\label{ex:cO}
Consider the case $k=l=n=3$. The minimal realisation is a partially truncated American football with $k=3$, which we 
can also think of as a painted tetrahedron, in which one face is black and the others are white. Another, non-minimal,
realisation is obtained by painting the faces of an octahedron so that two opposite faces are black, and the remaining 
six faces are white. This is not a branched cover of the painted tetrahedron.
\end{ex}
Note that in this example the minimal realisation has two kinds of vertices: one at which all faces meeting there are white, 
and three at which one black and two white faces meet. For the non-minimal realisation described above one has one 
black and three white faces meeting at every vertex. It is a general feature of the generalised football patterns with $n\geq 3$
that the combinatorial definition of the pattern does not imply any control over the local structure at a vertex. For $n=2$
we do have such control, leading to Theorem~\ref{t:cover}.

\bigskip

\bibliographystyle{amsplain}

\bigskip

\end{document}